\documentclass[12pt]{article}
\usepackage[english]{babel}
\usepackage[style=numeric,backend=biber]{biblatex}
\usepackage{indentfirst}
\usepackage{tikz}
\usetikzlibrary{calc}
\usepackage{amsmath,amssymb,amsthm}
\usepackage{cmap}
\usepackage{comment}
\usepackage{xcolor}
\usepackage{enumitem}

\usepackage[ruled,noline,titlenumbered,linesnumbered,noend,norelsize]{algorithm2e}
\DontPrintSemicolon
\SetNlSty{footnotesize}{}{}
\IncMargin{0.5em}
\SetNlSkip{1em}
\SetKwIF{If}{ElseIf}{Else}{if}{then}{else if}{else}{endif}
\SetKw{Return}{return}
\SetKw{Pick}{pick}
\SetNlSkip{2em}
\newenvironment{algorithm-hbox}{\hbadness=10000\begin{algorithm}}{\end{algorithm}}

\newenvironment{proof1}{\paragraph{Proof of Theorem 1.}}{\hfill$\square$}
\newenvironment{proof2}{\paragraph{Proof of Theorem 2.}}{\hfill$\square$}
\newenvironment{proof4}{\paragraph{Proof of Theorem 3. }}{\hfill$\square$}

\newtheorem{proposition}{Proposition}

\newcommand{\PP}{\mathbb{P}}

\newtheorem{theorem}{Theorem}

\usepackage[
    colorlinks=true,
    linkcolor=blue,
    citecolor=blue,
    urlcolor=blue
]{hyperref}

\begin{document}
\date{}

\title{Fractional hypergraph coloring}
\author{Margarita Akhmejanova \footnote{Matematiska institutionen, Uppsala universitet, Box 480, 751 06, Uppsala, Sweden. Email: mechmathrita@gmail.com}, Sean Longbrake \footnote{Emory University, Atlanta, GA 30322, USA. Email: sean.longbrake@emory.edu. }}

\maketitle

\begin{abstract}
We investigate proper $(a:b)$-fractional colorings of $n$-uniform hypergraphs, which generalize traditional integer colorings of graphs. Each vertex is assigned $b$ distinct colors from a set of $a$ colors, and an edge is properly colored if no single color is shared by all vertices of the edge.  A hypergraph is $(a:b)$-colorable if every edge is properly colored. We prove that for any $2\leq b\leq a-2$ and $a\leq n/(2\ln n)$, every $n$-uniform hypergraph $H$ with
$
|E(H)| \leq (ab^3)^{-1/2}\left(\frac{n}{\ln n}\right)^{1/2} \left(\frac{a}{b}\right)^{n-1}
$
is $(a:b)$-colorable. We also give a corresponding probabilistic upper bound and discuss specific cases.
\end{abstract}

\section{Introduction}

Graph coloring is a classical topic in combinatorics with numerous
applications in scheduling, coding theory, and resource allocation. In a
proper vertex coloring of a graph \(G\), each vertex receives one color and
adjacent vertices are required to receive different colors. Fractional
colorings of graphs can be introduced from two complementary points of view.

The first point of view comes from multicolorings. An \((a:b)\)-coloring of
a graph \(G\) assigns to each vertex a \(b\)-element subset of \([a]\), in such
a way that adjacent vertices receive disjoint subsets. The minimum such
integer \(a\) is called the \(b\)-fold chromatic number and is denoted by
\(\chi_b(G)\). The fractional chromatic number is then obtained by
normalizing these multicolorings:
\[
    \chi_f(G)=\inf_{b\geq 1}\frac{\chi_b(G)}{b}.
\]
This approach is closely related to the study of multicolorings of
graphs, notably by Stahl~\cite{Stahl1976}.

The second, equivalent, point of view comes from linear programming. Instead
of choosing a finite collection of color classes, one assigns nonnegative
weights to independent sets. If \(\mathcal I(G)\) denotes the family of
independent sets of \(G\), then
\[
    \chi_f(G)
    =
    \min \sum_{I\in \mathcal I(G)} w_I,
    \qquad
    \sum_{I\ni v} w_I\geq 1 \quad \text{for every } v\in V(G),
    \qquad
    w_I\geq 0.
\]
This linear programming interpretation connects fractional coloring with
fractional covers and was developed systematically in fractional graph
theory; see Lovász~\cite{lovasz1975ratio}, Grötschel, Lovász and
Schrijver~\cite{Grotschel1988}, and the monograph of Scheinerman and
Ullman~\cite{ScheinermanUllman1997}.

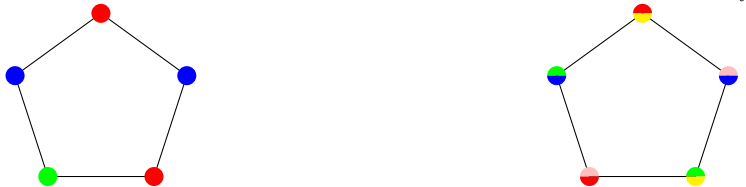
\begin{figure}[ht]
\caption{An example of \(3\)-proper and proper \((5:2)\)-fractional colorings of a pentagon graph \(G\). Here, \(\chi(G)=3\), while \(\chi_f(G)=2.5\).}
\label{picture:example}
\centering
\begin{tikzpicture}[line cap=round, line join=round, scale=1.2]

\coordinate (A1) at (90:1);
\coordinate (B1) at (18:1);
\coordinate (C1) at (-54:1);
\coordinate (D1) at (-126:1);
\coordinate (E1) at (162:1);

\draw (A1) -- (B1) -- (C1) -- (D1) -- (E1) -- cycle;

\fill[red]   (A1) circle (3pt);
\fill[blue]  (B1) circle (3pt);
\fill[red]   (C1) circle (3pt);
\fill[green] (D1) circle (3pt);
\fill[blue]  (E1) circle (3pt);


\begin{scope}[xshift=6cm]
  \coordinate (A2) at (90:1);
  \coordinate (B2) at (18:1);
  \coordinate (C2) at (-54:1);
  \coordinate (D2) at (-126:1);
  \coordinate (E2) at (162:1);

  \draw (A2) -- (B2) -- (C2) -- (D2) -- (E2) -- cycle;

  \fill[red]   (A2) -- ++(0:3pt) arc (0:180:3pt) -- cycle;
  \fill[yellow](A2) -- ++(180:3pt) arc (180:360:3pt) -- cycle;

  \fill[pink]  (B2) -- ++(0:3pt) arc (0:180:3pt) -- cycle;
  \fill[blue]  (B2) -- ++(180:3pt) arc (180:360:3pt) -- cycle;

  \fill[green] (C2) -- ++(0:3pt) arc (0:180:3pt) -- cycle;
  \fill[yellow](C2) -- ++(180:3pt) arc (180:360:3pt) -- cycle;

  \fill[pink]  (D2) -- ++(0:3pt) arc (0:180:3pt) -- cycle;
  \fill[red]   (D2) -- ++(180:3pt) arc (180:360:3pt) -- cycle;

  \fill[green] (E2) -- ++(0:3pt) arc (0:180:3pt) -- cycle;
  \fill[blue]  (E2) -- ++(180:3pt) arc (180:360:3pt) -- cycle;

\end{scope}
\end{tikzpicture}
\end{figure}

These interpretations also explain why fractional colorings naturally arise
in applications. In scheduling and timetabling problems, vertices may
represent jobs, and edges represent pairs of jobs that cannot be assigned to
the same time slot. If all jobs have unit processing time and each color
represents one time slot, then the chromatic number of the conflict graph is
the minimum number of time slots needed to complete all jobs without
conflicts. Fractional colorings refine this model by allowing jobs or
resources to be split among several time intervals. For instance, see
Figure~\ref{picture:example}, where a fractional coloring yields \(2.5\)
hours instead of \(3\). More general coloring-type scheduling frameworks
have also been studied; for example, Machacek introduced graph-like
scheduling problems defined by Boolean constraints of the form
\(x_i\neq x_j\), which generalize proper graph and hypergraph coloring
\cite{Machacek2025}. 

Fractional colorings also appear in coding theory,
where they are used to design efficient error-correcting codes and to
control conflicts between overlapping codewords; see
\cite{Korner1971, Malak2022}.

We now pass from graphs to hypergraphs. A hypergraph is a pair \(H=(V,E)\),
where \(V\) is the set of vertices and \(E\) is a family of subsets of \(V\),
called edges. A hypergraph is \(n\)-uniform if every edge contains exactly
\(n\) vertices. In the classical proper coloring of a hypergraph, each
vertex receives one color and no edge is allowed to be monochromatic.

In this paper we study an analogue of fractional graph coloring for uniform
hypergraphs. In a proper \((a:b)\)-coloring of a hypergraph, each vertex is
assigned \(b\) distinct colors from a set of \(a\) colors. An edge is
properly colored if no single color is shared by all vertices of the edge.
A hypergraph is \((a:b)\)-colorable if every edge is properly colored. When
all edges have size two, this definition agrees with the usual
\((a:b)\)-coloring of graphs.

This notion contains several classical coloring problems as special cases.
For \(b=1\), \((a:1)\)-colorability is the usual proper \(a\)-colorability
of hypergraphs. For \(b=a-1\), \((a:a-1)\)-colorability is equivalent to
panchromatic \(a\)-colorability. Recall that a panchromatic \(a\)-coloring
is a coloring in which each vertex receives a single color from \([a]\), and
every edge contains at least one vertex of each color. This equivalence can
be seen as follows. 
 In an $(a:a-1)$-coloring, each vertex receives exactly $a-1$ colors, thus missing precisely one color. If we assign each vertex its missing color, we establish a bijection between an $(a:a-1)$-coloring and a standard $a$-coloring where each vertex has exactly one color. A proper $(a:a-1)$-coloring ensures that every edge contains at least one vertex missing color $1$, at least one vertex missing color $2$, and so on, for all $a$ colors. Considering the missing colors, this directly transforms into a panchromatic $a$-coloring, where every edge contains all $a$ colors. Conversely, any panchromatic $a$-coloring can be transformed into a proper $(a:a-1)$-coloring by assigning each vertex all colors except its current one.

The classical Erd\H{o}s--Hajnal problem asks for the smallest number
\(m(n)\) of edges in an \(n\)-uniform hypergraph which is not properly
\(2\)-colorable. In 1963--1964, Erd\H{o}s~\cite{erdos1964, ErdLov} proved
that there exist absolute constants \(c,C>0\) such that
\[
    c2^n \leq m(n) \leq C2^n n^2 .
\]
Later, Radhakrishnan and Srinivasan~\cite{Radhakrishnan2000} improved the
lower bound by showing that there exists an absolute constant \(c'>0\) such
that
\[
    m(n) \geq c'2^n\left(\frac{n}{\ln n}\right)^{1/2}.
\]

The problem naturally extends to \(r>2\) colors. Let \(m(n,r)\) denote the
minimum number of edges in an \(n\)-uniform hypergraph which is not properly
\(r\)-colorable. For \(r\) small compared with \(n\), the best known bounds
have the form
\begin{equation}\label{bound:Cherk-Kozik}
c_1 \left(\frac{n}{\ln n}\right)^{\frac{r-1}{r}} r^{n-1}
\leq
m(n,r)
\leq
c_2 n^2 r^n\ln r,
\end{equation}
where \(c_1,c_2>0\) are absolute constants. The lower bound was proved by
Cherkashin and Kozik~\cite{CherkashinKozik2015}, while the upper bound
follows from a result of Erd\H{o}s; see the proof in Kostochka's
paper~\cite{Kostochka2002}. Bounds on \(m(n)\) for small values of \(n\), as
well as constructive examples, can be found, for example, in
\cite{cherkashin2023small, Aglave2016ImprovedBF, grill2024improved,
Radhakrishnan2021PropertyBT} and
\cite{gebauer2013construction, kozik2021improving}.

Kostochka~\cite{Kostochka2002} also considered the panchromatic analogue of
this problem. Namely, he studied the minimum number of edges in an
\(n\)-uniform hypergraph which does not admit a panchromatic \(r\)-coloring,
and denoted this number by \(p(n,r)\).  For more information and recent results, see \cite{Cherkashin2018, akhmejanova2022chain} and surveys \cite{raigorodskii2020extremal, raigorodskii2011erdHos}.

Fractional colorability of random hypergraphs has also been studied. A
random hypergraph \(H(n,k,p)\) has \(n\) vertices, and each \(k\)-element
subset is included as an edge independently with probability \(p\). In
\cite{zakharov2023, shabanov2024, kravtsov2019}, exponentially sharp
thresholds were obtained for the \((4:2)\)-, \((5:2)\)-, and
\((a:a-1)\)-colorability of \(H(n,k,p)\). These results are particularly
interesting because, for classical proper \(a\)-colorings with \(a>2\),
there remains a constant gap between the lower and upper bounds for the
threshold value of \(c\) when \(p=cn/\binom{n}{k}\).

\section{Our results}

\begin{theorem}\label{thm_1}
Assume that \(2\leq b\leq a-2\) and \(a\leq n/(2\ln n)\). Then every $n$-uniform hypergraph $H$ with
\begin{equation}\label{eq:thm_1}
|E(H)| \leq (ab^3)^{-1/2}\left(\frac{n}{\ln n}\right)^{1/2} \left(\frac{a}{b}\right)^{n-1}
\end{equation}
 is $(a:b)$-colorable.
\end{theorem}

Applying the Cherkashin--Kozik approach with ordered chains, as described
in \cite{CherkashinKozik2015}, one obtains the following bound: for
$\lfloor a/b\rfloor\geq 2$, every $n$-uniform hypergraph $H$ with
\begin{equation}\label{eq:Cherk-Kozik}
|E(H)| \leq
C_1
\left(\frac{n}{\ln n}\right)^{
\frac{\lfloor a/b\rfloor-1}{\lfloor a/b\rfloor}
}
\lfloor a/b\rfloor^{n-1}
\end{equation}
is properly $(a:b)$-colorable. Here $C_1>0$ is an absolute constant.

The bound (\ref{eq:Cherk-Kozik}) arises because the colors are divided into $k, k=\lfloor a / b \rfloor$ disjoint color groups 
$(1, 2, \ldots, b), (b+1, \ldots, 2b), \ldots, ((k-1)b+1, \ldots, kb)$, 
and vertices are then colored using these $k$-multicolors according to the Cherkashin–Kozik method. However, it is evident that (\ref{eq:Cherk-Kozik})
is weaker than (\ref{eq:thm_1}) when $a$ is not divisible by $b$ and $n$ is big enough. When $a$ is divisible by $b$ we receive lower bound from (\ref{bound:Cherk-Kozik}).

In Theorem \ref{thm_1} we may assume that $a - 2 \geq b$, as we explained before the case $a - 1 = b$ can be reduced to the problem of panchromatic coloring. In the paper \cite{akhmejanova2022chain} by the first author, together with Balogh and Shabanov, a bound for panchromatic colorings was proved. Hence, as a corollary, we immediately obtain the following result. 
\begin{proposition}
For $3 < a \leq \sqrt{\frac{n}{100 \ln n}}$, every $n$-uniform hypergraph $H$ with 
\begin{equation}
|E(H)| \leq \frac{1}{20a^2} \left(\frac{n}{\ln n}\right)^{\frac{a-1}{a}}\left(\frac{a}{a - 1}\right)^n
\end{equation}
is $(a:a-1)$-colorable.
\end{proposition}

We give an upper bound of the same general type as the classical upper bounds
for proper colorings and panchromatic colorings. Our argument follows
Erdős's probabilistic construction, but its extension to proper
$(a:b)$-colorings requires a more delicate estimate of the probability that a
random edge is monochromatic under a fixed coloring. The point is that a monochromatic random edge may have several common colors simultaneously.

\begin{theorem}\label{thm_2}
For any integers $a>b\ge 1$, there exists $n_0$, such as for all $n>n_0$ there is an $n$-uniform hypergraph with at most $\frac{e}{2}n^2\left(\frac{a}{b}\right)^n b(\ln\frac{a}{b}+1)$ edges which is not $(a:b)$-colorable.
\end{theorem}

If we suppose that $a/b\to\infty$, then Theorem \ref{thm_1} can be improved in the following way:

\begin{theorem}\label{thm_3}
Let
\(
    a'=\left\lfloor \frac{a(n-1)}{n}\right\rfloor.
\)
Suppose that \(b\leq a'\). If \(H=(V,E)\) is an \(n\)-uniform hypergraph with
\[
    |E(H)|\leq \left(\frac{a'}{b}\right)^n,
\]
then \(H\) is \((a:b)\)-colorable.
\end{theorem}

In Theorem \ref{thm_3} we use a reserve-color argument. We first color the vertices using only a subset of the available colors and keep the remaining
colors for a later recoloring step. A technical difficulty is that, in an
\((a:b)\)-coloring, an edge may have several colors common to all its
vertices. Thus, destroying a monochromatic edge may require deleting several
old colors. We overcome this by considering the corresponding bad pairs
\((A,\gamma)\), where \(\gamma\) is common to all vertices of an edge \(A\).
A Hall-type argument allows us to distribute these
bad pairs among vertices of the corresponding edges, after which the reserved
colors can be used without creating a new monochromatic edge.
\section{Proofs}
\begin{proof1}
Assume $H = (V, E)$ is a hypergraph with 
$|E(H)| \leq c\left(\frac{n}{\ln n}\right)^{1/2}\left(\frac{a}{b}\right)^n$, where constant $c = c(a, b)$ will be determined later.  To construct a proper $(a:b)$-fractional coloring of $H$, we proceed as follows.
First, assign to each vertex $v \in V$ independently of the other vertices a set of $b$ distinct colors chosen randomly from the set $\binom{[a]}{b}$. Next, assign to every vertex $v$ a random weight $\sigma(v)$ independently sampled from $U(0, 1)$, and reorder the vertices so that $\sigma(v_1) < \sigma(v_2) < \dots < \sigma(v_{|V|})$.  Finally, process vertices in the new order. Say, edge $A$ is monochromatic with color $\gamma\in[a]$, if all vertices of edge $A$ contain color $\gamma$ among the set of $b$ colors. If a vertex $v$ belongs to an initial monochromatic edge with color $\gamma_{\text{mod } a}$ and $\sigma(v) < p = \frac{1}{2} \frac{\ln \frac{n}{\ln n}}{n}$, replace the color $\gamma_{\text{mod } a}$ in $v$ with the next available color modulo $a$ so that all colors in $v$ remain distinct.  If $v$  is part of multiple initially  monochromatic edges with distinct colors, we randomly choose only one color to change at $v$.

We will show later that if the algorithm described above fails to produce a proper $(a:b)$-coloring for $H$, then one of the bad events $\mathcal{B}_1$, $\mathcal{B}_2$, $\mathcal{B}_3$, $\mathcal{B}_4$, $\mathcal{B}_5$ has occurred. Below, we describe these bad events.

\emph{Event} $\mathcal{B}_1$: there is a monochromatic edge $A$ with heavy vertices (vertices with weight more than $p$). 
$$\mathbb{P}(\mathcal{B}_1)\leq|E(H)|a\left(1-p\right)^n\left(\frac{\binom{a-1}{b-1}}{\binom{a}{b}}\right)^n\leq ac.$$
\emph{Event} $\mathcal{B}_2$: there exists a ``double" initially monochromatic edge $A$, meaning there are two distinct colors such that every vertex of $A$ contains both colors.
\begin{align*}
\mathbb{P}(\mathcal{B}_2)\leq
&a^2|E(H)|\left(\frac{\binom{a-2}{b-2}}{\binom{a}{b}}\right)^n\leq a^2c\sqrt{\frac{n}{\ln n}}\left(\frac{b-1}{a-1}\right)^n\\
&\leq\exp\left(\ln a^2c+\frac{\ln n}{2}-\frac{n(a-b)}{a-1}\right)\leq\frac{1}{n},
\end{align*}

when both conditions $a \leq \frac{n}{(2\ln n)}$ and $a^{3/2}c < 1$ are satisfied.

\emph{Event} $\mathcal{B}_3$: there is a light vertex $v$ (a vertex with weight less than $p$) that is part of two initially monochromatic edges $A$ and $B$ with distinct colors $\gamma_A$ and $\gamma_B$, respectively. 

We estimate the probability of this event by the union bound. There are at most
\(|E(H)|^2\) choices for the ordered pair \((A,B)\), at most \(a^2\)
choices for the colors \(\gamma_A,\gamma_B\), and at most \(k:=|A\cap B|\) choices
for the common vertex \(v\). We also apply that \(\sigma(v)<p\), and that common vertices contain both $\gamma_A$ and $\gamma_B$ colors and distinct vertices from $A$ and $B$ contain colors $\gamma_A$ and $\gamma_B,$ respectively. Hence,
\begin{gather*}
\mathbb{P}(\mathcal{B}_3) \leq 
|E(H)|^2 a^2 p \cdot
\max_k k
\left(\frac{\binom{a-2}{b-2}}{\binom{a}{b}}\right)^k
\left(\frac{\binom{a-1}{b-1}}{\binom{a}{b}}\right)^{2(n-k)} 
= \\
\frac{c^2 n}{\ln n}
\left(\frac{a}{b}\right)^{2n}
\frac{a^2}{2}\frac{\ln \frac{n}{\ln n}}{n}
\cdot
\max_k k
\left(\frac{b(b-1)}{a(a-1)}\right)^k 
\left(\frac{b}{a}\right)^{2(n-k)} 
\leq \\
\frac{c^2 a^2}{2}
\cdot
\max_k k
\left(\frac{1 - 1/b}{1 - 1/a}\right)^k 
=
\frac{c^2 a^2}{2}
\cdot
\max_k k 
\left(1 - \frac{a-b}{(a-1)b}\right)^k,
\quad x = \frac{a-b}{(a-1)b}, \\
\leq
\frac{c^2 a^2}{2}
\cdot
\max_k k e^{-kx} 
\leq
\frac{c^2 a^2}{2}
\cdot
\frac{1}{e x} 
\leq
\frac{a^3 c^2 b}{2 e (a-b)}.
\end{gather*}

We now consider the remaining possible type of failure: an edge $A$ becomes monochromatic with color
$\gamma_A$ during the recoloring procedure. Consider the recoloring step which
makes $A$ monochromatic with color $\gamma_A$, and let $v\in A$ be the vertex
recolored at this step. This recoloring is caused either by another initially
monochromatic edge $B\neq A$ with color $\gamma_B$, or by the edge $A$ itself. We treat these two
cases separately, since the corresponding probability estimates are different.

\emph{Event} $\mathcal{B}_4$: the first of the two cases above occurs, that is, the
vertex $v$ is recolored while processing an initially monochromatic edge
$B\neq A$  with color $\gamma_B$.
 
Let $\sigma(v)=x$ be the weight of $v$. 
Note several crucial facts by assuming that both events $\mathcal{B}_2$ and $\mathcal{B}_3$ do not occur:

\begin{enumerate}
    \item $v$ initially contains colors $\gamma_B, \gamma_B+1, \ldots, \gamma_A-1$ modulo $a$ and does not contain color $\gamma_A$.
    
    \item For all $u \in A \cap B$ with $u \neq v$, both colors $\gamma_A$ and $\gamma_B$ are initially present, and $\sigma(u) > \sigma(v)$. Indeed, $u$ must contain color $\gamma_B$ since the algorithm can only be initiated by initially monochromatic edges. From this and the assumption that $\mathcal{B}_1$ and $\mathcal{B}_2$ did not occur, we conclude that $\sigma(u) > \sigma(v)$ for all such $u$. Since $v$ is the last vertex recolored with $\gamma_A$ in $A$, all such $u$ must also contain color $\gamma_A$ in the initial coloring.
    
    \item Each $u \in B$, $u \neq v$, has weight at least $x$ and initially contains color $\gamma_B$.
    
    \item All $u \in A \setminus B$ with $\sigma(u) > x$ initially contain color $\gamma_A$ since $v$ is the last vertex recolored with $\gamma_A$ in $A$. All $u \in A \setminus B$ with $\sigma(u) < x$ initially contain either color $\gamma_A$ or $(\gamma_A - 1)$ modulo $a$, or possibly both.
\end{enumerate}

When both edges $A$ and $B$, vertex $v$, its weight $x$ and colors $\gamma_A, \gamma_B$ are fixed, we obtain the following estimate for the conditional probability:
\begin{gather}\label{eq:cond_prob}
\frac{\binom{a-t-1}{b-t}}{\binom{a}{b}}\left(\frac{\binom{a-2}{b-2}}{\binom{a}{b}}\right)^{k-1}(1-x)^{n-1}\left(\frac{\binom{a-1}{b-1}}{\binom{a}{b}}\right)^{n-k}\left(\frac{(1-x)\binom{a-1}{b-1}}{\binom{a}{b}}+\frac{2x\binom{a-1}{b-1}}{\binom{a}{b}}\right)^{n-k},
\end{gather}
where $k:=|A\cap B|$ and $ t:=(\gamma_A-\gamma_B)_{mod~a}.$ 

For a fixed \(\gamma_A\), the color \(\gamma_B\) can only be one of the
previous \(b\) colors in the cyclic order. Thus, 
we have \(1\leq t\leq b\).

By the Pascal triangle identity, also known as the hockey-stick identity,
we have
\[
    \sum_{t=1}^{b} \binom{a-t-1}{b-t}
    =
    \sum_{s=0}^{b-1} \binom{a-b+s-1}{s}
    =
    \binom{a-1}{b-1}.
\]
Therefore,
\[
    \sum_{t=1}^{b}
    \frac{\binom{a-t-1}{b-t}}{\binom{a}{b}}
    =
    \frac{\binom{a-1}{b-1}}{\binom{a}{b}}
    =
    \frac{b}{a}.
\]
To get the estimate for \(\mathbb{P}(\mathcal{B}_4)\), we integrate
\((\ref{eq:cond_prob})\) over the weight \(x\) and sum over all possible
choices of \(A,B,v,\gamma_A\), and \(t=1,\ldots,b\):
\begin{gather*}
\mathbb{P}\bigl(\mathcal{B}_4\setminus (\mathcal{B}_2\cup \mathcal{B}_3)\bigr)
\leq
|E|^2 a \cdot
\max_k k
\int_{0}^p 
\sum_{t=1}^b 
\frac{\binom{a-t-1}{b-t}}{\binom{a}{b}} 
\left(\frac{\binom{a-2}{b-2}}{\binom{a}{b}}\right)^{k-1}
(1-x)^{n-1}
\left(\frac{\binom{a-1}{b-1}}{\binom{a}{b}}\right)^{n-k}
\cdot \\
\cdot
\left(
\frac{(1-x)\binom{a-1}{b-1}}{\binom{a}{b}} 
+
\frac{2x\binom{a-1}{b-1}}{\binom{a}{b}}
\right)^{n-k} 
\, dx
\leq \\
b|E|^2
\int_{0}^p (1-x)^{n-1}(1+x)^{n-1} \, dx
\cdot
\max_k k
\left(\frac{b}{a}\right)^{2(n-k)}
\left(\frac{b(b-1)}{a(a-1)}\right)^{k-1}
\leq
\end{gather*}

\begin{gather*}
bp |E|^2 \cdot 
\max_k k \left(\frac{b}{a}\right)^{2(n-k)}
\left(\frac{b(b-1)}{a(a-1)}\right)^{k-1}
\leq \\
 \frac{b\ln \frac{n}{\ln n}}{2n}
\left(\frac{c^2 n}{\ln n}\right)
\left(\frac{a}{b}\right)^{2n}
\cdot
\max_k k \left(\frac{b}{a}\right)^{2(n-k)}
\left(\frac{b(b-1)}{a(a-1)}\right)^{k-1}
\leq \\
\frac{a^2 c^2}{2b}
\cdot
\max_k k
\left(\frac{1 - 1/b}{1 - 1/a}\right)^{k-1}
=
\frac{a^2 c^2}{2b}
\cdot
\max_k k
\left(1 - \frac{a-b}{(a-1)b}\right)^{k-1},
\quad
x=\frac{a-b}{(a-1)b}, \\
=
\frac{a^2 c^2}{2b}
\cdot
\max_k k(1-x)^{k-1}
\leq
\frac{a^2 c^2}{2b}
\cdot
\max_k k e^{-x(k-1)}
=
\frac{a^2 c^2e^x}{2b}
\cdot \max_k k e^{-kx}
\leq \\
\frac{a^2 c^2e^x}{2bex}
\leq
\frac{a^2 c^2}{2bx}
=
\frac{a^2 c^2(a-1)b}{2b(a-b)}
\leq
\frac{a^3c^2}{2(a-b)}.
\end{gather*}

We now turn to the case where the edge that becomes monochromatic is the
same edge that causes the recoloring step.

\emph{Event} $\mathcal{B}_5$: there exists an edge $A$ which was initially
monochromatic with some color $\gamma_B$, and after processing the edge
$A$ itself, the edge $A$ becomes monochromatic again, now with color
$\gamma_A$.

Assume that both events \(\mathcal{B}_2\) and \(\mathcal{B}_3\) do not occur.
Let \(v\in A\) be the vertex recolored while processing \(A\). Then \(v\) is
the first and only vertex of \(A\) that is recolored before \(A\) becomes
monochromatic with color \(\gamma_A\). Initially, all vertices of \(A\) had
color \(\gamma_B\), and therefore \(\gamma_A\neq\gamma_B\). There are at most \(|E(H)|\) choices for the edge \(A\), at most \(a^2\)
choices for the ordered pair of colors \((\gamma_A,\gamma_B)\), and at most
\(n\) choices for the vertex \(v\in A\). The vertex \(v\) must be light and
must initially contain the color \(\gamma_B\). Every other vertex
\(u\in A\setminus\{v\}\) must initially contain both colors \(\gamma_B\) and
\(\gamma_A\). Hence,
\begin{gather*}
\mathbb{P}(\mathcal{B}_5\setminus (\mathcal{B}_2\cup \mathcal{B}_3))\leq a^2|E(H)|np\left(\frac{\binom{a-2}{b-2}}{\binom{a}{b}}\right)^{n-1}\left(\frac{\binom{a-1}{b-1}}{\binom{a}{b}}\right)\leq a^2c\sqrt{n\ln n}\left(\frac{b-1}{a-1}\right)^{n-1}\\
\leq \exp\left(\ln a^2c+\frac{\ln n}{2} +\frac{\ln\ln n}{2}-\frac{(n-1)(a-b)}{a-1}\right)\leq\frac{1}{n},
\end{gather*}
provided that both conditions $a \leq \frac{n}{(2\ln n)}$ and $a^{3/2}c < 1$ are satisfied.

When the constant $c = (a^3b)^{-1/2}$, we have 
\[
\mathbb{P}(\cup_{i=1}^5\mathcal{B}_i) \leq ac +\frac{2}{n}+ \frac{a^3c^2b}{2e(a-b)} + \frac{a^3c^2}{2(a-b)} \leq \frac{1}{\sqrt{ab}}+\frac{2}{5}+\frac{1}{4e}+\frac{1}{8}< 1,\] since $a-b\geq 2$, $b\geq 2, a\geq 4$ and $n\geq 5$. 

So, with positive probability, the algorithm produces a coloring without these bad events. We claim that the algorithm produces a proper fractional coloring if none
of the events $\mathcal{B}_i$, $i\in[1,5]$, occurs. Indeed, the algorithm can fail only in two cases: either an initially
monochromatic edge remains monochromatic, which is covered by
$\mathcal{B}_1$, $\mathcal{B}_2$, and $\mathcal{B}_3$, or an edge becomes
monochromatic during the recoloring procedure, which is covered by
$\mathcal{B}_4$ and $\mathcal{B}_5$.
\end{proof1}

\begin{proof2}
We prove the theorem by adapting Erd\H{o}s's probabilistic argument for
proper \(2\)-colorings from \cite{erdos1964, ErdLov}. The main additional
difficulty is to control the probability that a random edge is bad under a
fixed \((a:b)\)-coloring. In contrast with the ordinary coloring case, such
an edge may be bad for several different colors simultaneously, and we
handle this overcounting by a careful inclusion--exclusion argument.

Let \(V\) be a set of size \(v\), where \(v\) will be chosen later. We first
allow repeated edges. Choose \(m\) independent uniformly random \(n\)-sets
\(S_1,\ldots,S_m\) from \(\binom{V}{n}\). We will show that, for a suitable
choice of \(m\), with positive probability the resulting \(n\)-uniform
multi-hypergraph has no proper \((a:b)\)-coloring. Deleting repeated edges
then gives an \(n\)-uniform hypergraph with at most \(m\) edges which is
still not \((a:b)\)-colorable.

Fix a map
\[
\chi:V\to \binom{[a]}{b}.
\]
We regard \(\chi\) as a possible \((a:b)\)-coloring of \(V\). Let \(S\) be a
uniformly random element of \(\binom{V}{n}\). We say that \(S\) is bad under
\(\chi\) if all vertices of \(S\) share at least one common color.

For a subset \(y\subseteq [a]\), define
\[
S_y=\{u\in V:y\subseteq \chi(u)\}.
\]
Thus \(S_y\) is the set of vertices whose color set contains all colors
from \(y\). 
Hence \(S\subseteq S_y\) if and only if all vertices of \(S\) have all
colors from \(y\) in common.

By the principle of inclusion and exclusion applied to the events
\(\{S\subseteq S_{\{x\}}\}\), \(x\in [a]\), we get
\[
\PP(S\text{ is bad under }\chi)
=
\sum_{i=1}^{b}(-1)^{i+1}
\frac{
\sum_{y\in\binom{[a]}{i}}\binom{|S_y|}{n}
}{
\binom{v}{n}
}.
\]
The sum stops at \(b\), since no vertex contains more than \(b\) colors.

Now fix a color \(x\in[a]\). For \(1\leq i\leq b\), put
\[
s_i^x=
\sum_{\substack{y\in\binom{[a]}{i}\\ x\in y}}
\binom{|S_y|}{n}.
\]
Then
\[
\sum_{y\in\binom{[a]}{i}}\binom{|S_y|}{n}
=
\frac{1}{i}\sum_{x\in[a]}s_i^x,
\]
because every \(i\)-set \(y\) contains exactly \(i\) colors. Therefore
\[
\PP(S\text{ is bad under }\chi)
=
\frac{1}{\binom{v}{n}}
\sum_{x\in[a]}\sum_{i=1}^b(-1)^{i+1}\frac{s_i^x}{i}.
\]

We claim that, for every fixed \(x\in[a]\),
\begin{equation}\label{eq:sx-bound}
\sum_{i=1}^b(-1)^{i+1}\frac{s_i^x}{i}
\geq
\frac{s_1^x}{b}.
\end{equation}
To prove this, we compare the contributions of each fixed
\(T\in\binom{V}{n}\) to both sides. For such \(T\), let
\[
R(T)=\bigcap_{u\in T}\chi(u)
\]
be the set of colors common to all vertices of \(T\), and let
\(
r(T)=|R(T)|.
\)
By the definition of \(S_y\), we have
\[
T\subseteq S_y
\quad\Longleftrightarrow\quad
y\subseteq R(T).
\]

If \(x\notin R(T)\), then \(T\) contributes zero to both sides of
\eqref{eq:sx-bound}. Suppose now that \(x\in R(T)\), and write
\(r=r(T)\). Then the contribution of \(T\) to \(s_i^x\) is
\[
\binom{r-1}{i-1},
\]
since the color \(x\) is already fixed and the remaining \(i-1\) colors are
chosen from the other \(r-1\) colors in \(R(T)\). Hence the contribution of
\(T\) to the left-hand side of \eqref{eq:sx-bound} is
\[
\sum_{i=1}^{r}(-1)^{i+1}\frac{1}{i}\binom{r-1}{i-1}
=
\frac1r\sum_{i=1}^{r}(-1)^{i+1}\binom{r}{i}
=
\frac1r.
\]
Since \(R(T)\subseteq \chi(u)\) for every \(u\in T\), and
\(|\chi(u)|=b\), we have \(r\leq b\). Thus \(1/r\geq 1/b\). On the other
hand, \(T\) contributes to \(s_1^x\) exactly once, and therefore its
contribution to \(s_1^x/b\) is \(1/b\). Thus the contribution of every
\(T\) to the left-hand side of \eqref{eq:sx-bound} is at least its
contribution to the right-hand side. Summing over all \(T\in\binom{V}{n}\)
proves \eqref{eq:sx-bound}.

It follows that
\[
\PP(S\text{ is bad under }\chi)
\geq
\frac{1}{b\binom{v}{n}}\sum_{x\in[a]}s_1^x.
\]
For \(x\in[a]\), let
\[
d_x=|S_{\{x\}}|.
\]
Thus \(d_x\) is the number of vertices whose color set contains \(x\), and
\(s_1^x=\binom{d_x}{n}\). Moreover,
\[
\sum_{x\in[a]}d_x=bv,
\]
since every vertex receives exactly \(b\) colors. Assume, for convenience,
that \(v\) is divisible by \(a\). Since \(k\mapsto \binom{k}{n}\) is convex
on the nonnegative integers, we obtain
\[
\sum_{x\in[a]}\binom{d_x}{n}
\geq
a\binom{bv/a}{n}.
\]
Hence, for every fixed \(\chi\),
\[
\PP(S\text{ is bad under }\chi)
\geq
p,
\qquad
p:=
\frac{a\binom{bv/a}{n}}{b\binom{v}{n}}.
\]

For a fixed coloring \(\chi\), the probability that \(\chi\) is proper for
all random edges \(S_1,\ldots,S_m\) is at most \((1-p)^m\), by independence.
Taking the union bound over all maps
\(\chi:V\to\binom{[a]}{b}\), we get
\[
\PP(\text{there exists a proper }(a:b)\text{-coloring})
\leq
\binom{a}{b}^v(1-p)^m.
\]
Using \(\binom{a}{b}\leq (ae/b)^b\) and \(1-p\leq e^{-p}\), this is at most
\[
\exp\left(bv\ln\frac{ae}{b}-pm\right).
\]
Choose
\[
m=
\left\lfloor
\frac{bv\ln\frac{ae}{b}}{p}
\right\rfloor+1.
\]
Then the last probability is strictly less than \(1\). Therefore there exists a
choice of \(S_1,\ldots,S_m\) such that the corresponding hypergraph is not
\((a:b)\)-colorable.

It remains to estimate \(m\). 
Rewrite $p$ as \[p = \frac{a \binom{\frac{bv}{a}}{n}}{b \binom{v}{n}} = \frac{a}{b} \prod_{i=0}^{n-1} \frac{\frac{bv}{a} - i}{v - i}=\frac{a}{b}\left( \frac{b}{a} \right)^n \prod_{i=0}^{n-1} \frac{1 - \frac{ai}{bv}}{1 - \frac{i}{v}} .\] 

Using \(\ln(1-z)=-z+O(z^2)\), we obtain, uniformly in \(i\),
\[
\ln\frac{1-\frac{ai}{bv}}{1-\frac{i}{v}}
=
\ln\left(1-\frac{ai}{bv}\right)
-
\ln\left(1-\frac{i}{v}\right)
=
-\frac{(a-b)i}{bv}
+
O\left(\frac{i^2}{v^2}\right).
\]
Hence
\[
\ln p
=
\ln\frac{a}{b}
+
n\ln\frac{b}{a}
-
\frac{a-b}{bv}\sum_{i=0}^{n-1}i
+
O\left(\sum_{i=0}^{n-1}\frac{i^2}{v^2}\right).
\]
If \(v=\Theta(n^2)\), the error term is \(o(1)\). Therefore
\[
p
\ge
(1-o(1))
\frac{a}{b}\left(\frac{b}{a}\right)^n
\exp\left(
-\frac{(a-b)n^2}{2bv}
\right).
\]

Finally, for \(v=\Theta(n^2)\), we conclude that
\[m  \leq (1+o(1))\frac{b^2v}{a} \left(\frac{a}{b}\right)^n \exp\left(\frac{(a-b)n^2}{2bv}\right) \ln \frac{ae}{b}.\]
Substituting \(v\), divisible by \(a\), such that
\(
v=(1+o(1))\frac{(a-b)n^2}{2b}
\), we obtain the desired bound.
\end{proof2}

\begin{proof4}
We use a reserve-color argument, in the spirit of Alon~\cite{Alon1985}, 
adapted to \((a:b)\)-colorings. First, we color the vertices using only the 
colors from \([a']\), while the remaining \(r\) colors \(a'+1,\ldots,a\) are 
kept for the recoloring step.

Assign to each vertex independently and uniformly at random a \(b\)-element
subset \(C(v)\subseteq [a']\). Let \(X\) is the number of pairs \((A,\gamma)\), where
\(A\in E(H)\), \(\gamma\in [a']\), and every vertex of \(A\) contains
\(\gamma\). We call such pairs bad pairs. By linearity of expectation,
\[
\mathbb{E}X
=
|E(H)|a'
\left(\frac{\binom{a'-1}{b-1}}{\binom{a'}{b}}\right)^n
=
|E(H)|a'\left(\frac{b}{a'}\right)^n
\leq a'.
\]
Hence, there exists an initial \((a':b)\)-coloring for which \(X\leq a'\).
Fix such a coloring.

We now assign each bad pair \((A,\gamma)\) to one vertex of the edge \(A\),
so that no vertex receives more than \(r\) assigned pairs. To prove that this
is possible, we construct an auxiliary bipartite graph \(G_{\mathrm{aux}}\).

The left part of \(G_{\mathrm{aux}}\) is the set of all bad pairs
\((A,\gamma)\). The right part consists of \(r\) copies of each vertex
\(v\in V(H)\), denoted by
\(
    v^{(1)},v^{(2)},\ldots,v^{(r)}.
\)
A bad pair \((A,\gamma)\) is adjacent to all copies of all vertices of \(A\);
that is, \((A,\gamma)\) is adjacent to \(v^{(j)}\) whenever \(v\in A\) and
\(1\leq j\leq r\).

We claim that \(G_{\mathrm{aux}}\) contains a matching which covers all bad
pairs. By Hall's theorem, it is enough to check that for every set \(S\) of bad
pairs,
\[
    |N(S)|\geq |S|,
\]
where \(N(S)\) denotes the set of neighbors of \(S\) in the auxiliary graph
\(G_{\mathrm{aux}}\). If \(S\neq\varnothing\), then \(S\) contains a pair \((A,\gamma)\), and hence
\(N(S)\) contains all \(r\) copies of every vertex of \(A\). Therefore
\[
    |N(S)|\geq nr.
\]
On the other hand,
\[
    |S|\leq X\leq a'\leq (n-1)r<nr.
\]

Thus Hall's condition is satisfied, and the desired matching exists.

Using this matching, we assign each bad pair \((A,\gamma)\) to the original
vertex \(v\in A\) whose copy \(v^{(j)}\) is matched to \((A,\gamma)\). Since
each vertex has only \(r\) copies, no vertex receives more than \(r\)
assigned bad pairs.

For every bad pair \((A,\gamma)\), delete the color \(\gamma\) from the
vertex assigned to this pair. If the same color is assigned several times to
the same vertex, it is deleted only once. After this step, no color from
\([a']\) is common to all vertices of any edge. Moreover, no vertex loses
more than \(\min\{r, b\}\) distinct colors.
It remains to replace the deleted colors by reserved colors. Let \(T\) be
the total number of distinct deleted color occurrences. Then
\[
    T\leq X\leq a'\leq (n-1)r.
\]
We choose the replacement colors so that each vertex receives distinct
reserved colors and each reserved color is used on at most \(n-1\) vertices.
Indeed, list the deleted color occurrences vertex by vertex, so that the occurrences
belonging to the same vertex form one consecutive block. Assign the reserved
colors \(\rho_1,\ldots,\rho_r\) cyclically along this list. Since each vertex
lost at most \(\min\{b,r\}\leq r\) colors, the colors assigned to the same
vertex are distinct. Moreover, under this cyclic assignment, each reserved color is used at most \(\lceil T/r\rceil\) times. Since \(T\leq (n-1)r\), each reserved color is used at most \(n-1\) times. After the replacement, no edge can be
monochromatic in a reserved color, since each reserved color is used on at
most \(n-1\) vertices in total.
\end{proof4}

\section{Acknowledgments}

The research work of M. Akhmejanova is supported by the Wenner-Gren Foundations.
 
\printbibliography

@article{erdos1964,
    author = {P. Erd\H{o}s},
    title = {On a combinatorial problem. II},
    journal = {Acta Math. Acad. Sci. Hungar.},
    volume = {15},
    number = {3-4},
    pages = {445--447},
    year = {1964}
}

@article{Radhakrishnan2000,
  author       = {J. Radhakrishnan and A. Srinivasan},
  title        = {Improved bounds and algorithms for hypergraph two-coloring},
  journal      = {Random Structures \& Algorithms},
  volume       = {16},
  number       = {1},
  year         = {2000},
  pages        = {4--32},
  doi          = {},
}

@article{Kostochka2002,
  author       = {A. Kostochka},
  title        = {On a theorem of Erd{\H{o}}s, Rubin, and Taylor on choosability of complete bipartite graphs},
  journal      = {Electronic Journal of Combinatorics},
  volume       = {9},
  number       = {1},
  year         = {2002},
  pages        = {1--4},
}

@incollection{ErdLov,
  author    = {P. Erd\H{o}s and L. Lov\'asz},
  title     = {Problems and results on 3-chromatic hypergraphs and some related questions},
  booktitle = {Infinite and Finite Sets},
  series    = {Colloquia Mathematica Societatis J\'anos Bolyai},
  volume    = {10},
  publisher = {North Holland},
  address   = {Amsterdam},
  year      = {1973},
  pages     = {609--627}
}

@article{CherkashinKozik2015,
  author    = {D. Cherkashin and J. Kozik},
  title     = {A note on random greedy coloring of uniform hypergraphs},
  journal   = {Random Structures and Algorithms},
  volume    = {47},
  number    = {3},
  year      = {2015},
  pages     = {407--413},
  doi       = {},
  publisher = {Wiley},
}

@article{Machacek2025,
  author  = {Machacek, John},
  title   = {Graph-Like Scheduling Problems and {Property B}},
  journal = {Studia Scientiarum Mathematicarum Hungarica},
  volume  = {62},
  number  = {1},
  pages   = {7--20},
  year    = {2025}
}

@article{Cherkashin2018,
  author  = {Cherkashin, D. D.},
  title   = {Panchromatic colorings of uniform hypergraphs},
  journal = {Journal of Combinatorial Theory, Series B},
  year    = {2018}
}

@article{Stahl1976,
  author  = {Stahl, Saul},
  title   = {n-tuple colorings and associated graphs},
  journal = {Journal of Combinatorial Theory, Series B},
  volume  = {20},
  number  = {2},
  pages   = {185--203},
  year    = {1976}
}

@book{ScheinermanUllman1997,
  author    = {Scheinerman, Edward R. and Ullman, Daniel H.},
  title     = {Fractional Graph Theory: A Rational Approach to the Theory of Graphs},
  publisher = {John Wiley \& Sons},
  address   = {New York},
  year      = {1997}
}

@article{Alon1985,
  author       = {N. Alon},
  title        = {Hypergraphs with high chromatic number},
  journal      = {Graphs and Combinatorics},
  volume       = {1},
  number       = {1},
  year         = {1985},
  pages        = {387--389},
  %doi          = {10.1007/BF02582913}
}

@article{akhmejanova2022chain,
  author       = {M. Akhmejanova and J. Balogh and D. Shabanov},
  title        = {Chain method for panchromatic colorings of hypergraphs},
  journal      = {Discrete Applied Mathematics},
  volume       = {321},
  pages        = {72--81},
  year         = {2022},
}

@article{zakharov2023,
  author       = {Zakharov, P. A. and Shabanov, D. A.},
  title        = {Fractional colourings of random hypergraphs},
  journal      = {Russian Mathematical Surveys},
  year         = {2023},
  volume       = {78},
  number       = {6},
  pages        = {183--184},
  %language     = {russian},
  note         = {(in Russian)},
  %doi          = {10.1070/RM10052} % Add DOI if available
}

@article{shabanov2024,
  author       = {D. Shabanov and T. Shaikheeva},
  title        = {Bounds for threshold probabilities for fractional
colorability properties of random hypergraphs},
  journal      = {Trudy MFTI},
  year         = {2024},
  volume       = {16},
  number       = {3},
  pages        = {81--90},
  note         = {(in Russian)},
  %language     = {russian},
}

@article{kravtsov2019,
  author       = {Kravtsov, D.  and Krokhmal, N. and Shabanov, D. },
  title        = {Panchromatic 3-colorings of random hypergraphs},
  journal      = {European Journal of Combinatorics},
  year         = {2019},
  volume       = {78},
  pages        = {28--43},
  % Add DOI if available
  keywords     = {random hypergraphs, panchromatic colorings, combinatorics}
}

@inproceedings{Malak2022,
  author    = {Malak, D.},
  title     = {Fractional Graph Coloring for Functional Compression with Side Information},
  booktitle = {Proceedings of the 2022 IEEE Information Theory Workshop (ITW)},
  year      = {2022},
  % Add any other fields if you have them, for example:
  % doi       = {...},
  % url       = {...},
  % pages     = {...},
}

@inproceedings{Korner1971,
  title={Coding of an information source having ambiguous alphabet and the entropy of graphs.},
  author={Korner, J.},
  booktitle={6th Prague conference on Information Theory, etc.},
  pages={411--425},
  year={1971},
  organization={Academia, Prague}
}

@article{cherkashin2023small,
  title={On small non-uniform hypergraphs without property B},
  author={Cherkashin, Danila},
  journal={Mathematics and Education in Mathematics},
  volume={52},
  pages={71--75},
  year={2023}
}

@article{Aglave2016ImprovedBF,
  title={Improved bounds for uniform hypergraphs without property B},
  author={Sachin, A. and V. Amarnath and Saswata Shannigrahi and Shwetank Singh},
  journal={Australas. J Comb.},
  year={2016},
  volume={76},
  pages={73-86},
}

@article{grill2024improved,
  title={Improved Lower Bounds for Property B},
  author={Grill, K. and Linzmayer, D.},
  journal={arXiv preprint arXiv:2403.05674},
  year={2024}
}

@inproceedings{Radhakrishnan2021PropertyBT,
  title={Property B: Two-Coloring Non-Uniform Hypergraphs},
  author={Radhakrishnan, J. and Srinivasan, S.},
  booktitle={Foundations of Software Technology and Theoretical Computer Science},
  year={2021},
}

@article{gebauer2013construction,
  title={On the construction of 3-chromatic hypergraphs with few edges},
  author={Gebauer, H.},
  journal={Journal of Combinatorial Theory, Series A},
  volume={120},
  number={7},
  pages={1483--1490},
  year={2013},
  publisher={Elsevier}
}

@inproceedings{kozik2021improving,
  title={Improving Gebauer’s Construction of 3-Chromatic Hypergraphs with Few Edges},
  author={Kozik, Jakub},
  booktitle={48th International Colloquium on Automata, Languages, and Programming (ICALP 2021)},
  year={2021},
  organization={Schloss-Dagstuhl-Leibniz Zentrum f{\"u}r Informatik}
}

@article{raigorodskii2020extremal,
  title={Extremal problems in hypergraph colourings},
  author={Raigorodskii, A. and Cherkashin, D.},
  journal={Russian Mathematical Surveys},
  volume={75},
  number={1},
  pages={89},
  year={2020},
  publisher={IOP Publishing}
}

@article{raigorodskii2011erdHos,
  title={The Erd{\H{o}}s-Hajnal problem of hypergraph colouring, its generalizations, and related problems},
  author={Raigorodskii, A. and Shabanov, D.},
  journal={Russian Mathematical Surveys},
  volume={66},
  number={5},
  pages={933},
  year={2011},
  publisher={IOP Publishing}
}

@book{Grotschel1988,
  author    = {M. Grötschel and L. Lovász and A. Schrijver},
  title     = {Geometric Algorithms and Combinatorial Optimization},
  publisher = {Springer},
  year      = {1988},
  address   = {Berlin},
  language  = {English}
}

@article{lovasz1975ratio,
  title={On the ratio of optimal integral and fractional covers},
  author={Lov{\'a}sz, L{\'a}szl{\'o}},
  journal={Discrete mathematics},
  volume={13},
  number={4},
  pages={383--390},
  year={1975},
  publisher={Elsevier}
}
\end{document}